\documentclass[12pt,reqno]{amsart}

\newtheorem{theorem}{Theorem}[section]

\theoremstyle{definition}   
\newtheorem{definition}{Definition}

\theoremstyle{remark}

\numberwithin{equation}{section}

\usepackage{times}
\usepackage{enumerate}
\usepackage{graphicx,adjustbox}
\usepackage{hyperref}
\usepackage{amsmath,amssymb}

\usepackage{amscd}
\usepackage{graphicx}
\usepackage[all]{xy}

\title[Numerical Semigroups Generated by Concatenation]
{Unboundedness of the first and the last Betti numbers of Numerical Semigroups Generated by Concatenation}
\author{
Ranjana Mehta
\and
Joydip Saha
\and
Indranath Sengupta
}
\date{}

\address{\small \rm  Department of Mathematics,
SRM University - AP, Amaravati, Andhra Pradesh 522502, India}
\email{ranjana.m@srmap.edu.in}

\address{\small \rm  Discipline of Mathematics, ISI Kolkata, Kolkata, 
West Bengal 700108, India.} 
\email{saha.joydip56@gmail.com}
\thanks{The second author thanks NBHM, Government of India for the post-doctoral fellow 
position at ISI Kolkata.}

\address{\small \rm  Discipline of Mathematics, IIT Gandhinagar, Palaj, Gandhinagar, 
Gujarat 382355, INDIA.}
\email{indranathsg@iitgn.ac.in}
\thanks{The third author is the corresponding author.}

\date{}



\allowdisplaybreaks

\begin{document}

\begin{abstract}
We show that the minimal number of generators and the Cohen-Macaulay type of a family of numerical semigroups generated by concatenation of arithmetic sequences is unbounded. 
\end{abstract}
\maketitle

\section{introduction}
A \textit{numerical semigroup} $\Gamma$ is a subset of the set of nonnegative 
integers $\mathbb{N}$, closed under addition, contains zero and generates 
$\mathbb{Z}$ as a group. We refer to \cite{rgs} for basic facts on numerical 
semigroups. In the paper \cite{mss3}, the authors have introduced the notion of concatenation of two arithmetic sequences to define a new family of numerical 
semigroups. This definition was largely inspired by the family of numerical 
semigroups in embedding dimension $4$, defined by Bresinsky in \cite{bre}. The 
authors have also proved in an earlier work \cite{mss2} that all the Betti 
numbers of the Bresinsky curves are unbounded. The whole motivation comes 
from the old question, whether every family of affine curves, in a fixed 
embedding dimension, has an upper bound on the minimal number of equations 
defining them ideal theoretically. The question was answered in negative 
for affine curves parametrised by monomials in \cite{bre}, and for 
algebroid space curves in \cite{moh}, \cite{mss3}. First example was probably 
given by F.S. Macaulay; we refer to \cite{abhyankar} for a detailed 
discussion on Macaylay's examples. However, all the examples found in the 
literature are in embedding dimension $4$ or lower. We did not find any 
example in embedding dimension $5$ or higher or for an arbitrary embedding 
dimension. This search led us to define the notion of concatenation of two numerical semigroups and we 
could indeed define a family in arbitrary embedding dimension in \cite{mss3}. 
In the same paper, we have conjectured that this family of numerical 
semigroups, in an arbitrary embedding dimension, defines affine monomial curves with 
arbitrarily large first Betti number (the minimal number of generators for the 
defining ideal) and the last Betti number (the Cohen-Macaulay type). In this 
article we show that, in embedding dimensions $4$ and $5$, our conjecture 
holds good.
\medskip

\section{Numerical semigroups and Concatenation}
Let us quickly discuss some basic facts on numerical semigroups and 
concatenation. Let $\Gamma$ be a numerical semigroup. It is true that 
(see \cite{rgs}) the set $\mathbb{N}\setminus \Gamma$ is 
finite and that the semigroup $\Gamma$ has a unique minimal system of generators 
$n_{0} < n_{1} < \cdots < n_{p}$. The greatest integer not belonging to $\Gamma$ 
is called the \textit{Frobenius number} of $\Gamma$, denoted by $F(\Gamma)$. The integers 
$n_{0}$ and $p + 1$ are known as the \textit{multiplicity} and the 
\textit{embedding dimension} of the semigroup $\Gamma$, usually 
denoted by $m(\Gamma)$ and $e(\Gamma)$ respectively. The 
\textit{Ap\'{e}ry set} of $\Gamma$ with respect to a non-zero $a\in \Gamma$ is 
defined to be the set $\rm{Ap}(\Gamma,a)=\{s\in \Gamma\mid s-a\notin \Gamma\}$. 
Given 
integers $n_{0} < n_{1} < \cdots < n_{p}$; the map 
$\nu : k[x_{0}, \ldots, x_{p}]\longrightarrow k[t]$ defined as 
$\nu(x_{i}) = t^{n_{i}}$, $0\leq i\leq p$, defines a parametrization 
for an affine monomial curve; the ideal $\ker(\nu)=\mathfrak{p}$ is called the 
defining ideal of the monomial curve defined by the parametrization 
$\nu(x_{i}) = t^{n_{i}}$, $0\leq i\leq p$. The defining ideal 
$\mathfrak{p}$ is a graded ideal with respect to the weighted gradation 
and therefore any two minimal generating sets of $\mathfrak{p}$ have the 
same cardinality. Similarly, by an abuse of notation, one can define a 
semigroup homomorphism $\nu: \mathbb{N}^{p+1} \rightarrow \mathbb{N}$ as \, 
$\nu((a_{0}, \ldots , a_{p})) = a_{0}n_{0} + a_{1}n_{1} + \cdots + a_{p}n_{p}$. 
Let $\sigma$ denote the kernel of congruence of the map $\nu$. It is known 
that $\sigma$ is finitely generated. The minimal number of generators of the 
ideal $\mathfrak{p}$, i.e., cardinality of a minimal generating set of 
$\mathfrak{p}$ is the same as the minimal cardinality of a system of generators 
of $\sigma$. 
\medskip

Let $\Gamma$ be a numerical semigroup, we say that $x\in\mathbb{Z}$ is a 
\textit{pseudo-Frobenius number} if $x\notin \Gamma$ and $x+s\in \Gamma$ 
for all $s\in \Gamma\setminus \{0\}$. We denote by $\mathbf{PF}(\Gamma)$ 
the set of  pseudo-Frobenius numbers of $\Gamma$. The cardinality of 
$\mathbf{PF}(\Gamma)$ is denoted by $t(\Gamma)$ and we call it the 
\textit{Cohen-Macaulay type} or simply the \textit{type} of $\Gamma$. 
Let $a,b\in \mathbb{Z}$. We define $\leq_{\Gamma}$ as $a\leq_{\Gamma} b$ 
if $b-a\in \Gamma$. This order relation defines a poset structure on 
$\mathbb{Z}$. It can be proved (see Proposition 8 in \cite{ap}) that 
$$\mathbf{PF}(\Gamma) = 
\{w-a\mid w\in \,\mathrm{Maximals}_{\leq_{\Gamma}}Ap(\Gamma, a)\}.$$ 

Let $e\geq 4$. Let us consider the string of positive integers 
in arithmetic progression: 
$a<a+d<a+2d<\ldots<a+(n-1)d<b<b+d<\ldots<b+(m-1)d$, where $m,n\in \mathbb{N}$, $m+n=e$  
and $\gcd(a,d)=1$. Note that $a<a+d<a+2d<\ldots<a+(n-1)d$ and 
$b<b+d<\ldots<b+(m-1)d$ are both arithmetic sequences with the 
same common difference $d$. We further assume that this sequence minimally generates 
the numerical semigroup $\Gamma=\langle a, a+d, a+2d,\ldots, a+(n-1)d, b, b+d,\ldots, b+(m-1)d\rangle$. 
Then, $\Gamma$ is called the \textit{numerical semigroup generated 
by concatenation of two arithmetic sequences with the same common difference $d$}. 
\medskip

Let $e\geq 4$, $n\geq 5$ and $q\geq 0$. We have defined in \cite{mss3} the 
numerical semigroup $\mathfrak{S}_{(n,e,q)}$, generated by the 
integers $\{m_{0},\ldots,m_{e-1}\}$, where 
$m_{i}:=n^2+(e-2)n+q+i$, for $0\leq i\leq e-3 $ and $m_{e-2}:=n^2+(e-1)n+q+(e-3)$, 
$m_{e-1}:=n^2+(e-1)n+q+(e-2)$. This is formed by concatenation of two arithmetic 
sequences with common difference $1$. Let  
$\mathcal{Q}_{(n,e,q)}\subset k[x_{0},\ldots,x_{e-1}]$ be the defining ideal of 
$\mathfrak{S}_{(n,e,q)}$. We have called $\mathfrak{S}_{(n,e,q)}$ an \textit{unbounded 
concatenation} and have conjectured the following:
\medskip

\noindent\textbf{Conjecture 3.3 \cite{mss3}.}\begin{enumerate}[(i)]
\item $\mu(\mathcal{Q}_{(n,e,e-4)})\geq n+2$;
\item the set $\{\mu(\mathcal{Q}_{(n,e,q)})\mid n\geq 5,e\geq 4, q\geq 0\}$ is unbounded above. 
\end{enumerate}
\medskip

\noindent We have also proved that $\mu(\mathcal{Q}_{(n,4,0)})= 2(n+1)$ in Theorem 3.5 in \cite{mss3}. 
This settles part (i) of the conjecture for the special case of $e=4$ and in deed 
justifies the naming \textit{unbounded concatenation}. 
\medskip

In this paper, we calculate the pseudo-Frobenius set 
for the case $e=4, q=0$ and we prove that the pseudo-Frobenius number or the 
Cohen-Macaulay type of the numerical 
semigroup ring is also unbounded. Note that the Cohen-Macaulay type is also the 
last Betti number. We then prove part (i) of the above conjecture for 
the case $e=5$ and this indeed gives us the desired class of monomial curves in embedding 
dimension $5$. We also compute the pseudo-Frobenius set and prove that the 
the pseudo-Frobenius number or the Cohen-Macaulay type in unbounded for $e=5$ 
as well. In light of these results we now revise the above conjecture and state 
the modified conjecture as follows:
\medskip

\noindent\textbf{Unbounded Concatenation Conjecture.}\begin{enumerate}[(i)]
\item $\mu(\mathcal{Q}_{(n,e,e-4)})\geq n+2$;
\item the Cohen-Macaulay type of $R/\mathcal{Q}_{(n,e,e-4)}$ is unbounded above;
\item the set $\{\mu(\mathcal{Q}_{(n,e,q)})\mid n\geq 5,e\geq 4, q\geq 0\}$ is unbounded above. 
\end{enumerate}
\bigskip

\section{Ap\'{e}ry set and the Pseudo-Frobenius set of $\mathfrak{S}_{(n,4)}$}
This section is devoted to the study of $\mathfrak{S}_{(n,e,q)}$, for $e=4$ and 
$q=e-4=0$. We will be writing $\mathfrak{S}_{(n,e)}$ instead of $\mathfrak{S}_{(n,e,0)}$,  
for simplicity of notation. Let $e\geq 4, i\geq 2$, $n=i(e-3)+(e-1)$ and 
$\mathfrak{S}_{(n,e,0)}=\langle m_{0},\ldots,m_{e-1}\rangle$, where 
\begin{align*}
m_{j} & = n^{2}+(e-2)n+(e-4+j),\quad 0\leq j\leq  e-3\\
m_{e-2} & = n^{2}+(e-1)n+(2e-7),\\
m_{e-1} & = n^{2}+(e-1)n+(2e-6).
\end{align*}

\begin{theorem}
Let $e\geq 4, i\geq 2$ and  $n=i(e-3)+(e-1)$, then the numerical semigroup $\mathfrak{S}_{(n,e)}$, is minimally generated by $\{m_{0},\ldots,m_{e-1}\}$ 
\end{theorem}
\proof See lemma 3.1 in \cite{mss3}\qed
\medskip

\noindent We note that $\mathfrak{S}_{(n,4)}=\langle m_{0},\ldots, m_{3}\rangle$, where $m_{0}=n^{2}+2n, m_{1}=n^{2}+2n+1, m_{2}=n^{2}+3n+1, m_{3}=n^{2}+3n+2$.
\medskip

\begin{theorem}\label{Apery 4}
The \textit{Ap\'{e}ry set} of $\mathfrak{S}_{(n,4)}$ with respect to  $m_{0}=n^{2}+2n$ is
$\rm{Ap}(\mathfrak{S}_{(n,4)},m_{0})=\displaystyle\cup_{i=1}^{5} A_{i}\cup \{0\}$. Where
\begin{align*}
A_{1}&=\{rm_{1}\mid 1\leq r\leq n\}\\
A_{2}&=\{rm_{2}\mid 1\leq r\leq n\}\\
A_{3}&=\{rm_{3}\mid 1\leq r\leq n-1\}\\
A_{4}&=\{rm_{1}+sm_{3}\mid 1\leq r\leq n-1,1\leq s\leq n-r\}\\
A_{5}&=\{rm_{2}+sm_{3}\mid 1\leq r\leq n-1,1\leq s\leq n-r\}.
\end{align*}
\end{theorem}

\proof \textbf{Case 1.} We show that $A_{1}\subset \rm{Ap}(\mathfrak{S}_{(n,4)},m_{0})$. We proceed by induction on $r$. For $r=1$, $m_{1}\in \rm{Ap}(\mathfrak{S}_{(n,4)},m_{0})$, because $m_{1}$ is an element of minimal generating set of $\mathfrak{S}_{(n,4)} $ hence it is true. Suppose $r>1$, we have $rm_{1}\equiv r(\rm{mod}\, m_{0})$, so there exists $s_{r}\in \rm{Ap}(\mathfrak{S}_{(n,4)},m_{0})$ such that $ rm_{1}=km_{0}+s_{r}$, $k\geq 0$ . Let $s_{r}=c_{1r}m_{1}+c_{2r}m_{2}+c_{3r}m_{3}$, then obviously $0\leq c_{1r}\leq r$. If 
$ c_{1r}=r $,  then $c_{2r}=c_{3r}=0$ and we are done. If $0\leq c_{1r}< r$, then $$ (r-c_{1r})m_{1}=km_{0}+c_{2r}m_{2}+c_{3r}m_{3} .$$
If $c_{1r}\neq 0$ then by induction, $(r-c_{1r})m_{1}\in \rm{Ap}(\mathfrak{S}_{(n,4)},m_{0})$ therefore $k=0$ and $rm_{1}=s_{r}\in \rm{Ap}(\mathfrak{S}_{(n,4)},m_{0})$.
We assume $c_{1r}=0$, which implies $s_{r}=c_{2r}m_{2}+c_{3r}m_{3}$. Hence 
\begin{align*}
rm_{1} &= km_{0}+c_{2r}m_{2}+c_{3r}m_{3},\quad  k\geq 0\\
(r-c_{2r}-c_{3r})m_{1} &= km_{0}+c_{2r}n+c_{3r}(n+1).
\end{align*}
As R.H.S. of the above equation is positive we have $c_{2r}+c_{3r}<r\leq n$. We have 
$$(r-c_{2r}-c_{3r}-k)m_{0}+(r-c_{2r}-c_{3r})=(c_{2r}+c_{3r})n+c_{3r}.$$
If $(r-c_{2r}-c_{3r}-k)=0 $, then from the above equation $(r-c_{2r}-c_{3r})=(c_{2r}+c_{3r})n+c_{3r} $ which gives a contradiction as $r<n$. Therefore $(r-c_{2r}-c_{3r}-k)\geq 1 $,implies $(r-c_{2r}-c_{3r}-k)m_{0}+(r-c_{2r}-c_{3r})\geq m_{1} $ hence we have $(c_{2r}+c_{3r})n+c_{3r}\geq m_{1}$. But $(c_{2r}+c_{3r})n+c_{3r}\leq rn+c_{3r}\leq n^{2}+n$ and $m_{1}=n^{2}+2n+1$. So we get a contradiction.
\medskip

\noindent\textbf{Case 2.} We want to show that $A_{2}\subset \rm{Ap}(\mathfrak{S}_{(n,4)},m_{0})$.
We have $m_{2}\in  \rm{Ap}(\mathfrak{S}_{(n,4)},m_{0})$. Let $r>1$, there exists $s_{r}\in  \rm{Ap}(\mathfrak{S}_{(n,4)},m_{0})$ such that $rm_{2}\equiv s_{r}(\rm{mod}\,m_{0})$. Therefore $rm_{2}=km_{0}+c_{1r}m_{1}+c_{3r}m_{3}$ $k\geq 0$ (by induction coefficient of $m_{2}$ say $c_{2r}$ is zero). Therefore $rm_{2}=km_{0}+c_{1r}m_{1}+c_{3r}(m_{2}+1)$, hence $(r-c_{3r})m_{2}= km_{0}+c_{1r}m_{1}+c_{3r}$. Since R.H.S. is positive we have $c_{3r}<r$. Again
\begin{align*}
r(m_{1}+n) & = km_{0}+c_{1r}m_{1}+c_{3r}(m_{1}+n+1)\\
(r-c_{1r}-c_{3r})m_{1} &= km_{0}+c_{3r}(n+1)-rn,\quad 1\leq r\leq n.
\end{align*}
If $k=0$ then $rm_{2}=s_{r}$ and we are done. We assume $k>0$. As R.H.S. is positive we get $c_{1r}+c_{3r}<r$. Now
\begin{align*}
& r(m_{0}+n+1) = km_{0}+c_{1r}(m_{0}+1)+c_{3r}(m_{0}+n+2)\\
\rm{implies}\quad &(r-k-c_{1r}-c_{3r})m_{0}+r(n+1)=c_{1r}+c_{3r}(n+2).
\end{align*}
If $r-k-c_{1r}-c_{3r}=0$ then $r(n+1)=(c_{1r}+c_{3r})+c_{3r}(n+1)<r+(r-1)(n+1)=r(n+1)+(r-n-1)$, $1\leq r\leq n$, which gives a contradiction.
\smallskip

If $r-k-c_{1r}-c_{3r}\neq 0$, then 
\begin{align*}
& (r-k-c_{1r}-c_{3r})m_{0}+r(n+1)\geq m_{0}+r(n+1)\\
\rm{implies}\quad & c_{1r}+c_{3r}(n+2)\geq m_{0}+r(n+1)\geq n^{2}+3n+1. 
\end{align*} 
But $c_{1r}+c_{3r}<r\leq n$, therefore $c_{1r}+c_{3r}(n+2)\leq n^{2}+2n-1$. Which gives a contradiction.
\medskip

\noindent\textbf{Case 3.} We wish to show that $A_{3}\subset \rm{Ap}(\mathfrak{S}_{(n,4)},m_{0})$. We note that for $r\geq 1$, $(r+1)m_{3}-m_{0}\notin \mathfrak{S}_{(n,4)}$ implies $rm_{3}-m_{0}\notin \mathfrak{S}_{(n,4)}$. Therefore it is enough to show that $(n-1)m_{3}-m_{0}\notin \mathfrak{S}_{(n,4)}$. Suppose 
\begin{align*}
(n-1)(m_{0}+n+2)-m_{0} & =c_{0}m_{0}+c_{1}m_{1}+c_{2}m_{2}+c_{3}m_{3}\\
 & =c_{0}m_{0}+c_{1}(m_{0}+1)+c_{2}(m_{0}+n+1)+c_{3}(m_{0}+n+2).
\end{align*} 
Therefore 
\begin{align*}
[(n-2)-(c_{0}+c_{1}+c_{2}+c_{3})]m_{0}+(n^{2}+n-2) &=c_{1}+c_{2}(n+1)+c_{3}(n+2).
\end{align*}
 As R.H.S. of the above equation is positive $(c_{0}+c_{1}+c_{3}+c_{4})\leq n-2.$ We have
 $$[(n-2)-(c_{0}+c_{1}+c_{2}+c_{3})]m_{0} =(c_{1}+c_{2}+2c_{3}+2)+n(c_{2}+c_{3}-n-1)$$ and $m_{0}=n(n+2)$. If $[(n-2)-(c_{0}+c_{1}+c_{2}+c_{3})]=0$ then we have $n(c_{2}+c_{3}-n-1)=-(c_{1}+c_{2}+2c_{3}+2)$, a contradiction, therefore $[(n-2)-(c_{0}+c_{1}+c_{2}+c_{3})]\neq 0$ and $n\mid (c_{1}+c_{2}+2c_{3}+2)$. Let $ c_{1}+c_{2}+2c_{3}+2=nk$ for some $k\geq 0$. If $k\geq 2$ then $c_{1}+c_{2}+2c_{3}+2\geq 2n$, on the other hand $c_{1}+c_{2}+c_{3}\leq n-2$ and $c_{3}\leq n-2$ implies $c_{1}+c_{2}+2c_{3}+2 \leq 2n-2<2n$, a contradiction. Therefore $k=0,1$. If $k=1$, then R.H.S. is $n(c_{2}+c_{3}-n)<0$, a contradiction. If $k=0$ then $c_{1}+c_{2}+2c_{3}+2=0$ again we get a contradiction.
\medskip
 
\noindent\textbf{Case 4.} We will show that $A_{4}\subset \rm{Ap}(\mathfrak{S}_{(n,4)},m_{0})$. We fix $r$, $1\leq r\leq n-1$. We need to show that  $rm_{1}+(n-r)m_{3}-m_{0}\notin \mathfrak{S}_{(n,4)}$. Suppose $rm_{1}+(n-r)m_{3}-m_{0}=\displaystyle\sum_{i=0}^{3}c_{i}m_{i}.$ After simplifying we get,
$$[n-(c_{0}+c_{1}+c_{2}+c_{3})]m_{0} =(c_{1}+c_{2}+2c_{3}+r)+n(c_{2}+c_{3}+r)$$ As R.H.S. is positive, we have $\displaystyle\sum_{i=0}^{3}c_{i}\leq n-1$ and $n\mid c_{1}+c_{2}+2c_{3}+r$. Let $c_{1}+c_{2}+2c_{3}+r=nk$, for some $k\geq 0$. Now $c_{1}+c_{2}+c_{3}\leq n-1, r\leq n-1, c_{3}\leq n-1$ implies $c_{1}+c_{2}+2c_{3}+r<3n $, hence $k\leq 2$. As $r>0$, we have $k\neq 0$. 
\smallskip

Suppose $k=1$ then $c_{1}+c_{2}+2c_{3}+r=n$ and R.H.S. is $n(c_{2}+c_{3}+r+1)$. Therefore $n+2\mid (c_{2}+c_{3}+r+1)$. Let $(c_{2}+c_{3}+r+1)=(n+2)\ell$ where $\ell>0$ (as $r+1>0$). We have $c_{2}+c_{3}\leq n-1$ and $r\leq n-1$, so $(c_{2}+c_{3}+r+1)<2(n+2)$. Hence $\ell=1$, We have $c_{2}+c_{3}+r=n+1$ and $c_{1}+c_{2}+2c_{3}+r=n$. Which implies $c_{1}+c_{3}+1=0$ a contradiction.
\smallskip

If $k=2$, then $c_{1}+c_{2}+2c_{3}+r=2n$ and R.H.S. is $n(c_{2}+c_{3}+r+2)$. Thus $n+2\mid (c_{2}+c_{3}+r+2)$. By the same arguments we get $c_{2}+c_{3}+r=n$. Therefore $c_{1}+c_{3}=n$ gives a contradiction as $c_{1}+c_{3}\leq n-1$.
\medskip

\noindent\textbf{Case 5.} We will show that $A_{5}\subset \rm{Ap}(\mathfrak{S}_{(n,4)},m_{0})$. We fix $r$, $1\leq r\leq n-1$. We need to show that  $rm_{2}+(n-r)m_{3}-m_{0}\notin \mathfrak{S}_{(n,4)}$. Suppose $rm_{2}+(n-r)m_{3}-m_{0}=\displaystyle\sum_{i=0}^{3}c_{i}m_{i}.$ After simplifying we get,
$$[n-( \displaystyle\sum_{i=0}^{3}c_{i})]m_{0}=(c_{1}+c_{2}+2c_{3}+r)+(c_{2}+c_{3})n.$$ As R.H.S. is positive we have $(\displaystyle\sum_{i=0}^{3}c_{i})\leq n-1$. Again we have $c_{1}+c_{2}+2c_{3}+r=kn$ with $k>0$. Now $c_{1}+c_{2}+c_{3}\leq n-1$, $c_{3}\leq n-1$ and $r\leq n-1$ implies that $c_{1}+c_{2}+2c_{3}+r<3n$. Therefore $k=1,2$.
\smallskip

If $k=1$, $c_{1}+c_{2}+2c_{3}+r=n$ and R.H.S. is $n(c_{2}+c_{3}+1)$. Now $n+2\mid c_{2}+c_{3}+1$. Let $(c_{2}+c_{3}+1)=(n+2)\ell$, where $\ell\neq 0$. Which gives a contradiction as $c_{2}+c_{3}\leq n-1$.
\smallskip

If $k=2$, $c_{1}+c_{2}+2c_{3}+r=2n$ and R.H.S. is $n(c_{2}+c_{3}+2)$. Now $n+2\mid c_{2}+c_{3}+2$. Let $(c_{2}+c_{3}+2)=(n+2)\ell$, where $\ell\neq 0$. Again a contradiction as $c_{2}+c_{3}\leq n-1$. \qed
\medskip

\begin{theorem}\label{pf4}
The pseudo-Frobenius set of the numerical semigroup $\mathfrak{S}_{(n,4)}$ is $PF(\mathfrak{S}_{(n,4)})=P_{1}\cup P_{2}\cup P_{3}$, where
\begin{align*}
P_{1}&=\{(n-1)m_{0}+n\};\\
P_{2}&=\{(n-1)m_{0}+n+k(n+1)\mid  1\leq k\leq (n-1)\};\\
P_{3}&=\{(n-1)m_{0}+n+(n-1)(n+1)+t\mid 1\leq t \leq n\}.
\end{align*}
\end{theorem}

\proof Let $P^{'}_{j}=P_{j}+m_{0}$, $1\leq j\leq 3$, therefore  
\begin{align*}
P^{'}_{1}&= \{nm_{1}\};\\
P^{'}_{2}&= \{km_{1}+(n-k)m_{3},1\leq k\leq n-1\};\\
P^{'}_{3}&= \{km_{2}+(n-k)m_{3},0\leq k\leq n-1\}.
\end{align*}
We want to show the following two conditions:
\smallskip
\begin{enumerate}
\item[(i)] For each $x\in Ap(\mathfrak{S}_{(n,4)},m_{0})\setminus \{P^{'}_{1}\cup P^{'}_{2} \cup P^{'}_{3}\}$, there exists  $ y\in \{P^{'}_{1}\cup P^{'}_{2} \cup P^{'}_{3}\}$ such that $y-x \in \mathfrak{S}_{(n,4)}$.

\item[(ii)] For any $y_{1},y_{2}\in \{P^{'}_{1}\cup P^{'}_{2} \cup P^{'}_{3}\}$, $y_{1}-y_{2} \notin \mathfrak{S}_{(n,4)}$. 
\end{enumerate}
\medskip

\noindent\textbf{Proof of (i):} Let $x=rm_{1}$, where $1\leq r\leq n-1$. Take $y=nm_{1}$ and $y-x=(n-r)m_{1}\in \mathfrak{S}_{(n,4)}$.
\medskip

Let $x=rm_{2}$, where $1\leq r\leq n-1$ then take $y=rm_{2}+(n-r)m_{3}\in\{P^{'}_{1}\cup P^{'}_{2} \cup P^{'}_{3} $ and $y-x=(n-r)m_{3}\in \mathfrak{S}_{(n,4)}$.
\medskip

For $x\in A_{3}\subset  Ap(\mathfrak{S}_{(n,4)},m_{0})\setminus \{P^{'}_{1}\cup P^{'}_{2} \cup P^{'}_{3}\}$, we take $y=m_{1}+(n-1)m_{3}$. Then we have 
$ y-x=m_{1}+(n-1-r)m_{3}\in \mathfrak{S}_{(n,4)}$ for $1\leq r\leq n-1$.
\medskip

 For $x\in  Ap(\mathfrak{S}_{(n,4)},m_{0})\setminus \{P^{'}_{1}\cup P^{'}_{2} \cup P^{'}_{3}\}$ and $x\in\{rm_{1}+sm_{3}|1\leq r \leq n-2, 1\leq s\leq n-r-1\}$ $x=rm_{1}+s m_{3}$ then we choose  $y=rm_{1}+(n-r) m_{3}$ so that we have $y-x=(n-r-s)m_{3}$ since $1\leq s\leq n-1-r$ therefore $(n-r-s)\geq 1$
\medskip

For $x\in  Ap(\mathfrak{S}_{(n,4)},m_{0})\setminus \{P^{'}_{1}\cup P^{'}_{2} \cup P^{'}_{3}\}$ and $x\in\{rm_{2}+s m_{3}|1\leq r \leq n-2, 1\leq s \leq n-r-1\}$ $x=rm_{2}+s m_{3}$ then we choose  $y=rm_{2}+(n-r) m_{3}$ so that we have $y-x=(n-r-s)m_{2}$ since $1\leq s\leq n-1-r$ therefore $(n-r-s)\geq 1$.
\medskip

\noindent\textbf{Proof of (ii):} It is easy to check for any $y_{1},y_{2}\in P^{'}_{j}$, $1\leq j\leq 3$, $y_{1}-y_{2}\notin\mathfrak{S}_{(n,4)}$.
\medskip

Let $y_{1}= nm_{1}$ and $y_{2}\in P^{'}_{2}$, then
\begin{align*}
y_{2}-y_{1}&=km_{1}+(n-k)m_{3}-nm_{1}\\
&=(k-n)m_{1}+(n-k)(m_{1}+n+1)\\
&=n^{2}+n-k(n+1)
\end{align*}
Since $1\leq k \leq n-1$ therefore $y_{2}-y_{1}=n^{2}+n-k(n+1)< m_{0}$ therefore $y_{1}-y_{2}\notin\mathfrak{S}_{(n,4)}$.
\medskip

Let $y_{1}= nm_{1}$ and $y_{2}\in P^{'}_{3}$, then
\begin{align*}
y_{2}-y_{1}&=k m_{2}+(n-k)m_{3}-nm_{1}\\
&=(k-n)m_{1}+kn+(n-k)(m_{1}+n+1)\\
&=n^2+n-k
\end{align*}
Since $0\leq k \leq n-1$ therefore $y_{2}-y_{1}=n^2+n-k< m_{0}$ therefore $y_{1}-y_{2}\notin\mathfrak{S}_{(n,4)}$.
\medskip

Let $y_{1}=\in P^{'}_{2}$ and $y_{2}\in P^{'}_{3}$, then
\begin{align*}
y_{2}-y_{1}&=k_{2}m_{2}+(n-k_{2})m_{3}-k_{1}m_{1}-(n-k_{1})m_{3}\\
&=(k_{2}-k_{1})m_{1}+k_{2}n+(k_{1}-k_{2})m_{3}\\
&=k_{2}n+(k_{1}-k_{2})(n+1)\\
&=k_{1}(n+1)
\end{align*}
Since $1\leq k_{1} \leq n-1$ therefore $y_{2}-y_{1}=k_{1}(n+1)< m_{0}$ therefore $y_{1}-y_{2}\notin\mathfrak{S}_{(n,4)}$.
 \qed
\medskip

\section{Ap\'{e}ry set and the Pseudo-Felonious set of $\mathfrak{S}_{(n,5)}$}
We now consider the numerical semigroup $\mathfrak{S}_{(n,5,0)}$. We 
write $\mathfrak{S}_{(n,5)}$ instead of $\mathfrak{S}_{(n,5,0)}$, 
for simplicity of notation. We first describe the Ap\'{e}ry 
set of $\mathfrak{S}_{(n,5)}$, which will help us write the 
pseudo-Frobenius set. We have proved that the pseudo-Frobenius 
number of the Cohen-Macaulay type is indeed unbounded. In the 
next section, we find a minimal generating set for the defining 
ideal and show that the minimal number of generators of the 
defining ideal is unbounded above, thereby proving parts (i) and (ii) of the Unbounded Concatenation Conjecture. 
\medskip

Let us recall that for $e\geq 4, i\geq 2$, $n=i(e-3)+(e-1)$, the numerical semigroup 
$\mathfrak{S}_{(n,e,0)}=\langle m_{0},\ldots,m_{e-1}\rangle$ is minimally 
generated by $m_{0},\ldots,m_{e-1}$, where 
\begin{align*}
m_{j} & = n^{2}+(e-2)n+(e-4+j),\quad 0\leq j\leq  e-3\\
m_{e-2} & = n^{2}+(e-1)n+(2e-7),\\
m_{e-1} & = n^{2}+(e-1)n+(2e-6).
\end{align*}

\noindent We note that $\mathfrak{S}_{(n,5)}=\langle m_{0},\ldots, m_{4}\rangle$, 
where $m_{0}=n^{2}+3n+1$, $m_{1}=n^{2}+3n+2$, $m_{2}=n^{2}+3n+3$, 
$m_{3}=n^{2}+4n+3$, $m_{4}= n^{2}+3n+4$. 
Let $\mathcal{Q}_{(n,5)}\subset k[x_{0},\ldots,x_{4}]$ be the defining ideal of 
$\mathfrak{S}_{(n,5)}$.
\medskip

\begin{theorem}\label{Apery 5}
The \textit{Ap\'{e}ry set} of $\mathfrak{S}_{(n,5)}$ with respect to  $m_{0}=n^{2}+3n+1$ is
$\rm{Ap}(\mathfrak{S}_{(n,5)},m_{0})=\displaystyle\cup_{i=1}^{11} A_{i}$. Where
\begin{align*}
A_{1}& =\{0,m_{1}\}\\
A_{2}&=\{rm_{2}\mid 1\leq r\leq \dfrac{n}{2}\}\\
A_{3}&=\{rm_{3}\mid 1\leq r\leq n\}\\
A_{4}&=\{rm_{4}\mid 1\leq r\leq n\}\\
A_{5}&=\{m_{1}+r m_{2}\mid 1\leq r\leq \dfrac{n}{2}\}\\
A_{6}&=\{m_{3}+r m_{2}\mid 1\leq r\leq \dfrac{n}{2}\}\\
A_{7}&=\{rm_{2}+2s m_{4}\mid 1\leq s\leq \dfrac{n}{2}-1,1\leq r \leq \dfrac{n}{2}-s \}\\
A_{8}&=\{rm_{2}+(2s-1)m_{4}\mid 1\leq s\leq \dfrac{n}{2},1\leq r \leq \dfrac{n}{2}+1-s \}\\
A_{9}&=\{rm_{3}+(n-k-r+1)m_{4}\mid 1\leq k\leq n-1,1\leq r \leq n-k \}\\
A_{10}&=\{rm_{2}+m_{3}+2s m_{4}\mid 1\leq s\leq \dfrac{n}{2}-1,1\leq r \leq \dfrac{n}{2}-s \}\\
A_{11}&=\{rm_{2}+m_{3}+(2s-1)m_{4}\mid 1\leq s\leq \dfrac{n}{2}-1,1\leq r \leq \dfrac{n}{2}+1-s \}.
\end{align*}
\end{theorem}
\proof

\noindent \textbf{Case 1.} Clearly $A_{1}\in \rm{Ap}(\mathfrak{S}_{(n,5)},m_{0})$
\medskip

\noindent \textbf{Case 2.} To show $A_{2}, A_{3}, A_{4}\in \rm{Ap}(\mathfrak{S}_{(n,5)},m_{0})$, we proceed similarly as in the cases 1, 2 and 3 of \ref{Apery 4}.
\medskip 

\noindent \textbf{Case 3.} To show $A_{5}\in \rm{Ap}(\mathfrak{S}_{(n,5)},m_{0})$, 
it is enough to show that $m_{1}+\dfrac{n}{2}m_{2}-m_{0}\notin \mathfrak{S}_{(n,5)}$. 
Let $m_{1}+\dfrac{n}{2}m_{2}-m_{0}=c_{0}m_{0}+c_{1}m_{1}+c_{2}m_{2}+c_{3}m_{3}+c_{4}m_{4}$. Converting $m_{i}, 1\leq i\leq 4$ in the term of $m_{0}$, we get the 
following equation 
$$(\dfrac{n}{2}-\sum\limits_{i=0}^{4}c_{i})m_{0}=c_{1}+2c_{2}+(n+2)c_{3}+(n+3)c_{4}-(n+1)$$
If $\sum\limits_{i=0}^{4}c_{i}<\dfrac{n}{2}$, then $(\dfrac{n}{2}-\sum\limits_{i=0}^{4}c_{i})m_{0}> m_{0}$, and\\[2mm]
$c_{1}+2c_{2}+(n+2)c_{3}+(n+3)c_{4}-(n+1)$\\[2mm]
$ < c_{1}+2c_{2}+(n+2)c_{3}+(n+3)c_{4}$\\[2mm]
$ \leq c_{1}+2c_{2}+ 2c_{3}+3c_{4}+\dfrac{n^{2}}{2}$\\[2mm]
$< \dfrac{n^{2}}{2}+\dfrac{3n}{2}< m_{0}$,\\[2mm]
which is a contradiction. 
\medskip

\noindent Again if $\sum\limits_{i=0}^{4}c_{i}> \dfrac{n}{2}$, then $(\dfrac{n}{2}-\sum\limits_{i=0}^{4}c_{i})m_{0}\leq -m_{0}$, and\\
 $c_{1}+2c_{2}+(n+2)c_{3}+(n+3)c_{4}-(n+1)\geq -(n+1)>-m_{0}$, which gives a contradiction.
\medskip

\noindent If $\sum\limits_{i=0}^{4}c_{i}= \dfrac{n}{2}$, then $c_{1}+2c_{2}+(n+2)c_{3}+(n+3)c_{4}=n+1$, which gives $(c_{0}+c_{1}+c_{2}+c_{3}+c_{4})+ c_{2}+(n+1)c_{3}+(n+2)c_{4}=n+1-c_{0}$, hence $c_{2}+(n+1)c_{3}+(n+2)c_{4}=1+\dfrac{n}{2}-c_{0}$\\
Proceeding in the similar way we get, $nc_{3}+(n+1)c_{4}=1-2c_{0}-c_{1}$, and  $c_{4}=1-\dfrac{n^2}{2}-(n+2)c_{0}-(n+1)c_{1}-nc_{2}$.\\
Since $n\geq 8$ and $c_{0}, c_{1}, c_{2}\geq 0$, the above equation gives $c_{4}<0$, which is not possible. 
\smallskip

\noindent Similarly we can show that $A_{6}\in \rm{Ap}(\mathfrak{S}_{(n,5)},m_{0})$.
\medskip

\noindent \textbf{Case 4.} To show $A_{7}\in \rm{Ap}(\mathfrak{S}_{(n,5)},m_{0})$. At first we fix $1\leq s\leq \dfrac{n}{2}-1$.\\ We need to show $(\dfrac{n}{2}-s)m_{2}+2s m_{4}-m_{0}\notin \mathfrak{S}_{(n,5)}$.\\
Let $(\dfrac{n}{2}-s)m_{2}+2s m_{4}-m_{0}=c_{0}m_{0}+c_{1}m_{1}+c_{2}m_{2}+c_{3}m_{3}+c_{4}m_{4}$, which gives,\\  
$(\dfrac{n}{2}+s-1-\sum\limits_{i=0}^{4}c_{i})m_{0}=c_{1}+2c_{2}+(n+2)c_{3}+(n+3)c_{4}-n(2s+1)-4s$\\
If $\dfrac{n}{2}+s-1>\sum\limits_{i=0}^{4}c_{i}$, then 
$(\dfrac{n}{2}+s-1-\sum\limits_{i=0}^{4}c_{i})m_{0}\geq m_{0}$ and\\[2mm]
$c_{1}+2c_{2}+(n+2)c_{3}+(n+3)c_{4}-n(2s+1)-4s$\\[2mm]
$\leq c_{1}+2c_{2}+2c_{3}+3c_{4}+(n+3)(\dfrac{n}{2}+s-1)-n(2s+1)-4s$\\[2mm]
$\leq 3(\dfrac{n}{2}+s-1)+\dfrac{n^2}{2}-\dfrac{n}{2}-ns-s-3$\\[2mm]
$\leq\dfrac{n^2}{2}+n+2s\leq \dfrac{n^2}{2}+2n<m_{0}$, 
which is a contradiction. 
\medskip

\noindent If $\dfrac{n}{2}+s-1<\sum\limits_{i=0}^{4}c_{i}$, then $(\dfrac{n}{2}+s-1-\sum\limits_{i=0}^{4}c_{i})m_{0}\leq -m_{0}$, and\\[2mm]
$c_{1}+2c_{2}+(n+2)c_{3}+(n+3)c_{4}-n(2s+1)-4s$\\[2mm]
$\geq -n(2s+1)-4s$\\[2mm]
$\geq -n^{2}-n+4 \quad (\text{substituting}\, s=\dfrac{n}{2}-1)$\\[2mm]
$ > -m_{0}$, \, which is a contradiction. 
\medskip

\noindent If $\dfrac{n}{2}+s-1=\sum\limits_{i=0}^{4}c_{i}$, then $c_{1}+2c_{2}+(n+2)c_{3}+(n+3)c_{4}=n(2s+1)+4s$. 
We have $(\dfrac{n}{2}-s)m_{2}+2s m_{4}=c_{0}^{'}m_{0}+c_{1}m_{1}+c_{2}m_{2}+c_{3}m_{3}+c_{4}m_{4}$, where $c_{0}^{'}=c_{0}+1\geq 1$. Now 
\begin{align*}
(\dfrac{n}{2}-s)m_{2}+2s m_{4}&=n^{2}(\dfrac{n}{2}+s)+(n+1)(\dfrac{3n}{2}+5s)\\
&=(n^{2}-1)(\dfrac{n}{2}+s)+(n+1)(\dfrac{3n}{2}+5s)+(\dfrac{n}{2}+s)\\
&=(n+1)[\dfrac{n^{2}}{2}+n+4s+ns]+(\dfrac{n}{2}+s)\equiv (\dfrac{n}{2}+s) mod (n+1)
\end{align*}
and\\[2mm]
$c'_{0}(m_{1}-1)+c_{1}m_{1}+c_{2}(m_{1}+1)+c_{3}(m_{1}+n+1)+c_{4}(m_{1}+n+2)$\\[2mm]
$(c_{2}+c_{4}-c'_{0})mod (n+1)$,\\[2mm]
which implies $(\dfrac{n}{2}+s)\equiv (c_{2}+c_{4}-c'_{0})mod (n+1)$. Substituting 
$c'_{0}+c_{1}+c_{2}+c_{3}+c_{4}= \dfrac{n}{2}+s$, we get 
$c'_{0}+c_{1}+c_{2}+c_{3}+c_{4}\equiv (c_{2}+c_{4}-c'_{0}) mod (n+1)$. Therefore, $n+1$ divides $2c'_{0}+c_{1}+c_{3}$, which implies $2c'_{0}+c_{1}+c_{3}=0$ (since $c_{2}+c_{4}-c_{0}^{'}> 0$). Hence $c'_{0}=0$, which is a contradiction to $c'_{0}\geq 1$.
\smallskip

\noindent Since expressions of $A_{8}, A_{9}, A_{10}, A_{11}$ are similar to the expression of $A_{7}$, therefore proofs follow similar steps.\qed
\medskip

\begin{theorem}
The set of all pseudo Frobenius numbers of the numerical semigroup $\mathfrak{S}_{(n,5)}$ is 
$$PF(\mathfrak{S}_{(n,5)})=P_{1}\cup P_{2}\cup P_{3}.$$ Where, 
\begin{align*}
P_{1}&=\{\dfrac{n}{2}m_{0}+(n+1)\}\\
P_{2}&=\{\dfrac{n}{2}m_{0}+(n+1)+k(m_{3}+(n+2))\mid 1\leq k\leq (\dfrac{n}{2}-1)\}\\
P_{3}&=\{\dfrac{n}{2}m_{0}+(n+1)+(\dfrac{n}{2}-1)(m_{3}+n+2)+(n+1+t)\mid 0\leq t\leq n+2\}.
\end{align*}
\end{theorem}
\proof \proof Similar as \ref{pf4}.\qed
\medskip

\section{minimal generating set for the defining ideal $\mathcal{Q}_{(n,5)}$}
\begin{theorem}\label{gastinger}
 Let $A = k[x_{1},\ldots,x_{n}]$ be a polynomial ring, $I\subset A$ the defining
ideal of a monomial curve defined by natural numbers $a_{1},\ldots,a_{n}$, 
whose greatest common divisor is $1$.  Let $J \subset I$ be a subideal. 
Then $J = I$ if and only if $\mathrm{dim}_{k} A/\langle J + (x_{i}) \rangle =a_{i}$
for some $i$. (Note that the above conditions are also equivalent to 
$\mathrm{dim}_{k} A/\langle J + (x_{i}) \rangle =a_{i}$ for any $i$.)
\end{theorem}  

\proof See \cite{g}.\qed
\medskip
 
\begin{theorem}
Let $e\geq 4, i\geq 2$ and  $n=i(e-3)+(e-1)$ then a minimal generating set of the defining ideal of $\mathfrak{S}_{(n,5)}$ consist of the following polynomials, 
\begin{itemize}
\item $f_{1}=x_{1}x_{3}-x_{0}x_{4}$.
\item $f_{2}=x_{2}x_{3}-x_{1}x_{4}$.
\item $g_{1}=x_{1}^2-x_{0}x_{2}$.
\item $g_{2}=x_{2}^{i+3}-x_{0}^{i+2}x_{3}$.
\item $\xi_{t}=x_{0}^tx_{1}^{n+2-t}-x_{3}^{t}x_{4}^{n+1-t},\quad 0\leq t \leq n+1$.
\item $\eta_{k}=x_{0}^{k+1}x_{3}^{n-2k-1}-x_{2}^{k+2}x_{4}^{n-2k-2},\quad 0\leq k \leq i$.
\item $l_{1}=x_{0}^{n+1}x_{1}-x_{2}x_{4}^n$.
\item $l_{2}=x_{0}^{n+2}-x_{2}x_{3}x_{4}^{n-1}$.
\end{itemize}
 \end{theorem}
 \proof
We consider the set $S=\{f_{1},f_{2},g_{1},g_{2}, \xi_{t}, \eta_{k}, l_{1}, l_{2}\mid0\leq t\leq n+1,\, 0\leq k\leq i,\}$. Let $J$ be the ideal generated by $S$, we consider the ideal $J+\langle x_{0}\rangle$. Then a generating set of $J+\langle x_{0}\rangle$ is 
\begin{itemize}
\item $q=x_{0}$
\item $\tilde{f_{1}}=x_{1}x_{3}$
\item $\tilde{f_{2}}=x_{2}x_{3}-x_{1}x_{4}$
\item $\tilde{g_{1}}=x_{1}^2$
\item $\tilde{g_{2}}=x_{2}^{i+3}$
\item $\tilde{\xi_{0}}=x_{4}^{n+1}$
\item $\tilde{\xi_{t}}=-x_{3}^{t}x_{4}^{n+1-t},\quad 1\leq t \leq n+1$
\item $\tilde{\eta_{k}}=-x_{2}^{k+2}x_{4}^{n-2k-2},\quad 0\leq k \leq i$
\item $\tilde{l_{1}}=-x_{2}x_{4}^n$
\item $\tilde{l_{2}}=-x_{2}x_{3}x_{4}^{n-1}$
\end{itemize}
We consider the lexicographic monomial order induced by $x_{0}\geq x_{1}\geq x_{2}\geq x_{3}\geq x_{4}$ on $k[x_{0},\ldots,x_{4}]$. Then a a standard basis of $J+\langle x_{0}\rangle$ w.r.t the given monomial order consist of following polynomials
 \begin{itemize}
\item $q=x_{0}$
\item $\tilde{f_{1}}=x_{1}x_{3}$
\item $\tilde{f_{2}}=x_{2}x_{3}-x_{1}x_{4}$
\item $\tilde{g_{1}}=x_{1}^2$
\item $\tilde{g_{2}}=x_{2}^{i+3}$
\item $\tilde{\xi_{0}}=x_{4}^{n+1}$
\item $\tilde{\xi_{t}}=x_{3}^{t}x_{4}^{n+1-t},\quad 1\leq t \leq n+1$
\item $\tilde{\eta_{k}}=x_{2}^{k+2}x_{4}^{n-2k-2},\quad 0\leq k \leq i$
\item $\tilde{l_{1}}=x_{2}x_{4}^n$
\item $\tilde{l_{2}}=x_{2}x_{3}x_{4}^{n-1}$
\item $h=x_{2}x_{3}^2$ 
\end{itemize}

Since all the generators of ideal $J+\langle x_{0}\rangle$ except $\tilde{f}_{2}$ are monomials therefore we calculate only all the S-polynomials with $\tilde{f}_{2}$ and they are as follows:
\begin{itemize}
\item $S(\tilde{f_{2}}, \tilde{f_{1}})=-x_{2}x_{3}^2=-h$
\item $S(\tilde{f_{2}}, \tilde{g_{1}})=-x_{1}x_{2}x_{3}=-x_{2}\cdot\tilde{f_{1}}$
\item It is clear that $\gcd(\rm{Lt}(\tilde{f_{2}}), \tilde{g_{2}})=1$
\item  $S(\tilde{f_{2}}, \tilde{\xi_{0}})=-x_{2}x_{3}x_{4}^{n}=-x_{2}\cdot \tilde{\xi_{1}}$
\item For $1\leq t\leq n$, we have  $S(\tilde{f_{2}}, \tilde{\xi_{t}})=x_{2}x_{3}^{t+1}=x_{3}^{t}\cdot\tilde{f_{2}}+x_{3}^{t-1}x_{4}\cdot\tilde{f_{1}}$
\item Since $\xi_{n+1}=x_{3}^{n+1}$, hence $\gcd(\rm{Lt}(\tilde{f_{2}}), \tilde{\xi_{n+1}})=1$
\item For $0\leq k\leq i$,  $S(\tilde{f_{2}}, \tilde{\eta_{k}})=-x_{2}^{k+3}x_{3}x_{4}^{n-2k-2}=-x_{2}^{k+3}\cdot\xi_{0}$
\item  $S(\tilde{f_{2}}, \tilde{l_{1}})=-x_{2}^{2}x_{3}x_{4}^{n-1}=x_{2}\cdot\tilde{l_{2}}$
\item $S(\tilde{f_{2}}, \tilde{l_{2}})=-x_{2}^{2}x_{3}^2x_{4}^{n-2}=x_{2}x_{3}x_{4}^{n-2}\cdot\tilde{f_{2}}-x_{1}\cdot\tilde{l_{2}}$
\end{itemize}
Therefore the set $$T=\{q,\tilde{f_{1}},\tilde{f_{2}},\tilde{g_{1}},\tilde{g_{2}},\tilde{\xi_{0}},\tilde{\xi_{t}},\tilde{\eta_{k}},\tilde{l_{1}},\tilde{l_{2}},h\mid 1\leq t\leq n+1, 0\leq k\leq i\}$$ forms a standard basis for the ideal $J+\langle x_{0}\rangle$. Hence, the leading ideal 
$\mathrm{lead}(J+\langle x_{0}\rangle)$  of 
$J+\langle x_{0}\rangle$, with respect to the given monomial order is generated by the 
following set, 
\begin{eqnarray*}
G & = & \{x_{0}, x_{1}x_{3}, x_{1}x_{4}, x_{1}^2, x_{2}^{i+3}, x_{4}^{n+1}, x_{2}x_{4}^{n}, 
x_{2}x_{3}x_{4}^{n-1}\} \cup\\
{} & {} &  \{x_{3}^{t}x_{4}^{n+1-t},x_{2}^{k+2}x_{4}^{n-2k-2} \mid 1\leq t \leq n+1,0\leq k \leq i \}.
\end{eqnarray*}
We need to show that $\dim_{k}\left(k[x_{0},\ldots,x_{4}]/J+\langle x_{0}\rangle\right)=m_{0}$. 
We list all monomials which are not divided by any element of $G$.
\begin{itemize}
\item $\{1,x_{1}\}$
\item $\{x_{2},x_{2}^2, \ldots,x_{2}^{i+2}\}$
\item $\{x_{3},x_{3}^{2},\ldots, x_{3}^n\}$
\item $\{x_{4},x_{4}^2,\ldots, x_{4}^n\}$
\item $\{x_{1}x_{2}, x_{1}x_{2}^2,\ldots, x_{1}x_{2}^{i+2}\}$
\item $\{x_{2}x_{3},x_{2}^2x_{3},\ldots,x_{2}^{i+2}x_{3}\}$
\item \begin{align*}
&\{x_{2}x_{4},\ldots, x_{2}x_{4}^{n-1}\}\\
&\{x_{2}^2x_{4},\ldots, x_{2}^2x_{4}^{n-3}\}\\
&\vdots\\
&\{x_{2}^{i+2}x_{4},\ldots, x_{2}^{i+2}x_{4}^{n-2i-3}\}\\
\end{align*}
\item \noindent\begin{align*}
&\{x_{3}x_{4},\ldots,x_{3}x_{4}^{n-1}\}\\
&\{x_{3}^{2}x_{4},\ldots, x_{3}^2x_{4}^{n-2}\}\\
&\vdots\\
&\{x_{3}^{n-1}x_{4}\}
\end{align*}
\item \noindent\begin{align*}
&\{x_{2}x_{3}x_{4},\ldots, x_{2}x_{3}x_{4}^{n-2}\}\\
&\{x_{2}^{2}x_{3}x_{4},\ldots, x_{2}^{2}x_{3}x_{4}^{n-3}\}\\
&\{x_{2}^{2}x_{3}x_{4},\ldots, x_{2}^{2}x_{3}x_{4}^{n-5}\}\\
&\vdots\\
&\{x_{2}^{i+2}x_{3}x_{4},\ldots, x_{2}^{i+2}x_{3}x_{4}^{n-2i-3}\}
\end{align*}
\end{itemize}
The cardinality of the set containing all the above elements given is \\
$2+2i+2n+(i+2)+(i+2)+(i+2)(n-i-2)+\dfrac{n(n-1)}{2}+(i+2)(n-i-2)-1$.\\
Using $(i+2)=\dfrac{n}{2}$ we get\\
$2(n+1)+\dfrac{3n}{2}+\dfrac{n^{2}}{2}-\dfrac{n^{2}}{4}+\dfrac{n^{2}}{2}- \dfrac{n}{2}+\dfrac{n^{2}}{2}-\dfrac{n^{2}}{4}-1=n^{2}+3n+1=m_{0}$.\\
Therefore, by Theorem \ref{gastinger}, the set $S$ is a generating set for the defining ideal 
of $\mathfrak{S}_{(n,5)}$.
\medskip

To prove minimality of the generating set, we prove that no element of this generating set can be expressed by other elements. 
\medskip

Let $\pi_{024}: k[x_{0},x_{1},x_{2},x_{3},x_{4}]\rightarrow k[x_{1},x_{3}]$ be defined as,\\
$\pi_{024}(x_{0})=0, \pi_{024}(x_{1})=x_{1}, \pi_{024}(x_{2})=0, \pi_{024}(x_{3})=x_{3},\pi_{024}(x_{4})=0$. 
Let
$$f_{1}= x_{1}x_{3}-x_{0}x_{4}= c_{f2} f_{2}+c_{g1}g_{1}+c_{g2}g_{2}+\sum\limits_{t=0}^{t=n+1}c_{\xi t}\xi_{t}+\sum\limits_{k=0}^{k=i}c_{\eta k}\eta_{k}+c_{l1} l_{1}+c_{l2}l_{2},$$
for $c_{f2}, c_{g1}, c_{g2}, c_{\xi t}, c_{\eta k}, c_{l1}, c_{l2}\in k[x_{0},x_{1},x_{2},x_{3},x_{4}]$. Applying $\pi_{024}$ to both the sides of the equation we get,  
$x_{1}x_{3}=c_{g1}(0,x_{1},0,x_{3},0)x_{1}^{2}$. By comparing the exponent of 
$x_{1}$ on both sides of the above equation we can conclude that the above equation 
is not possible.
\medskip

Let 
$\pi_{014}: k[x_{0},x_{1},x_{2},x_{3},x_{4}]\rightarrow k[x_{2},x_{3}]$ be defined as,\\
 $\pi_{014}(x_{0})=0, \pi_{014}(x_{1})=0, \pi_{014}(x_{2})=x_{2}, \pi_{014}(x_{3})=x_{3},\pi_{014}(x_{4})=0$. Let 
$$f_{2}= x_{2}x_{3}-x_{1}x_{4}= c_{f1} f_{1}+c_{g1}g_{1}+c_{g2}g_{2}+\sum\limits_{t=0}^{t=n+1}c_{\xi t}\xi_{t}+\sum\limits_{k=0}^{k=i}c_{\eta k}\eta_{k}+c_{l1} l_{1}+c_{l2}l_{2},$$
for $c_{f1}, c_{g1}, c_{g2}, c_{\xi t}, c_{\eta k}, c_{l1}, c_{l2}\in k[x_{0},x_{1},x_{2},x_{3},x_{4}]$. Applying $\pi_{014}$ to both the sides of the equation we get,  
$x_{2}x_{3}=c_{g2}(0,0,x_{2},x_{3},0)x_{2}^{i+3}$, where $i\geq 2$. 
By comparing the exponent of $x_{2}$ in both the sides of the above equation we can 
conclude that the above equation is not possible. 
\medskip

Let
$$g_{1}= x_{1}^{2}-x_{0}x_{2}= c_{f1} f_{1}+c_{f2}f_{2}+c_{g2}g_{2}+\sum\limits_{t=0}^{t=n+1}c_{\xi t}\xi_{t}+\sum\limits_{k=0}^{k=i}c_{\eta k}\eta_{k}+c_{l1} l_{1}+c_{l2}l_{2},$$
for $c_{f1}, c_{f2}, c_{g2}, c_{\xi t}, c_{\eta k}, c_{l1}, c_{l2}\in k[x_{0},x_{1},x_{2},x_{3},x_{4}]$. 
Applying $\pi_{024}$ to both the sides of the equation we get, $x_{1}^{2}= c_{f1}(0,x_{1},x_{2},x_{3},0) (x_{1}x_{3})$, where $i\geq 2$. By comparing the exponent of $x_{3}$ in both the sides of the 
above equation we can conclude that the above equation is not possible.
\medskip

Let
$$g_{2}= x_{2}^{i+3}-x_{0}^{i+2}x_{3}= c_{f1} f_{1}+c_{f2}f_{2}+c_{g1}g_{1}+\sum\limits_{t=0}^{t=n+1}c_{\xi t}\xi_{t}+\sum\limits_{k=0}^{k=i}c_{\eta k}\eta_{k}+c_{l1} l_{1}+c_{l2}l_{2},$$
for $c_{f2}, c_{g1}, c_{\xi t}, c_{\eta k}, c_{l1}, c_{l2}\in k[x_{0},x_{1},x_{2},x_{3},x_{4}]$. 
Applying $\pi_{014}$ to both the sides of the equation we get, 
$x_{2}^{i+3}=  c_{f2}(0,0,x_{2},x_{3},0) (x_{2}x_{3})$, where $i\geq 2$. 
By comparing the exponent of $x_{3}$ in both the sides of the above equation we can conclude 
that the above equation is not possible. 
\medskip

Let $\pi_{01}: k[x_{0},x_{1},x_{2},x_{3},x_{4}]\rightarrow k[x_{2},x_{3},x_{4}]$ be defined as,\\ 
$\pi_{01}(x_{0})=0, \pi_{01}(x_{1})=0, \pi_{01}(x_{2})=x_{2}, \pi_{01}(x_{3})=x_{3},\pi_{01}(x_{4})=x_{4}$.\\
Let 
$$\xi_{t}=x_{0}^tx_{1}^{n+2-t}-x_{3}^{t}x_{4}^{n+1-t}=c_{f1} f_{1}+c_{f2}f_{2}+c_{g1}g_{1}+c_{g2}g_{2}+\sum\limits_{k=0}^{k=i}c_{\eta k}\eta_{k}+c_{l1} l_{1}+c_{l2}l_{2}.$$
Applying $\pi_{01}$ to both the sides of the equation we get, 
\begin{eqnarray*}
-x_{3}^{t}x_{4}^{n+1-t} & = & c_{f2}(0,0,x_{2},x_{3},x_{4})(x_{2}x_{3})+c_{g2}(0,0,x_{2},x_{3},x_{4})x_{2}^{i+3}\\
{} & {} & -\sum\limits_{k=0}^{k=i}c_{\eta k}(0,0,x_{2},x_{3},x_{4}) x_{2}^{k+2}x_{4}^{n-2k-2}-c_{l1}(0,0,x_{2},x_{3},x_{4}) x_{2}x_{4}^{n}\\
{} & {} & -c_{l2}(0,0,x_{2},x_{3},x_{4})(x_{2}x_{3}x_{4}^{n-1}), \quad 0\leq t\leq n+1.
\end{eqnarray*}
From the above equations we observe that every term in the R.H.S contains $x_{2}$ but L.H.S does not have $x_{2}$, therefore the above equation is not possible.
\medskip

Let $\pi_{124}: k[x_{0}, x_{1}, x_{2}, x_{3}, x_{4}]\rightarrow k[x_{0}, x_{3}]$ be defined as,\\
$\pi_{124}(x_{0}) = x_{0}, \pi_{124}(x_{1}) = 0, \pi_{124}(x_{2}) = 0, \pi_{124}(x_{3}) = x_{3}, 
\pi_{124}(x_{4})=0$. Let 
$$\eta_{k}=x_{0}^{k+1}x_{3}^{n-2k-1}-x_{2}^{k+2}x_{4}^{n-2k-2}=c_{f1} f_{1}+c_{f2}f_{2}+c_{g1}g_{1}+c_{g2}g_{2}+\sum\limits_{t=0}^{t=n+1}c_{\xi t}\xi_{t}+c_{l1} l_{1}+c_{l2}l_{2}.$$ 
Applying $\pi_{124}$ to both the sides of the equation we get:
$$x_{0}^{k+1}x_{3}^{n-2k-1}= -c_{g2}(x_{0},0,0,x_{3},0)x_{0}^{i+2}x_{3}+c_{l2}(x_{0},0,0,x_{3},0)x_{0}^{n+2}, \quad 0\leq k\leq i.$$
Since $0\leq k\leq i$, $i\geq 2$ and $k+1 < i+2$, $n=4+2i$, the above equation is not possible.  
\bigskip

Let $\pi_{013}: k[x_{0},x_{1},x_{2},x_{3},x_{4}]\rightarrow k[x_{2},x_{4}]$ be defined as,\\
$\pi_{013}(x_{0}) = 0, \pi_{013}(x_{1}) = 0, \pi_{013}(x_{2}) = x_{2}, \pi_{013}(x_{3}) = 0, 
\pi_{013}(x_{4})=x_{4}$. Let 
$$l_{1}= x_{0}^{n+1}x_{1}-x_{2}x_{4}^n =c_{f1} f_{1}+c_{f2}f_{2}+c_{g1}g_{1}+c_{g2}g_{2}+\sum\limits_{t=0}^{t=n+1}c_{\eta k}\eta_{k}+\sum\limits_{t=0}^{t=n+1}c_{\xi t}\xi_{t}+c_{l2}l_{2}.$$
Applying $\pi_{013}$ to both the sides of the equation we get:\\
$$-x_{2}x_{4}^n= c_{g2}(0,0,x_{2},0,x_{4})x_{2}^{i+3}-\sum\limits_{k=0}^{k=i}c_{\eta k}(0,0,x_{2},0,x_{4})x_{2}^{k+2}x_{4}^{n-2k-2}.$$
Since $i\geq 2$, each term of R.H.S of the above equation contains $x_{2}^{2}$ whereas L.H.S containing only $x_{2}$. Hence the above equation is not possible.
\bigskip

Let $\pi_{234}: k[x_{0},x_{1},x_{2},x_{3},x_{4}]\rightarrow k[x_{0},x_{1}]$ be defined as,\\
$\pi_{234}(x_{0})=x_{0}, \pi_{234}(x_{1})=x_{1}, \pi_{234}(x_{2})=0, \pi_{234}(x_{3})=0, \pi_{234}(x_{4})=0$. Let
$$l_{2}= x_{0}^{n+2}-x_{2}x_{3}x_{4}^{n-1}=c_{f1} f_{1}+c_{f2}f_{2}+c_{g1}g_{1}+c_{g2}g_{2}+\sum\limits_{k=0}^{k=i}c_{\eta k}\eta_{k}+\sum\limits_{t=0}^{t=n+1}c_{\xi t}\xi_{t}+c_{l1}l_{1}.$$ 
Applying $\pi_{234}$ to both the sides of the equation we get:
$$x_{0}^{n+2}= c_{g1}(x_{0},x_{1},0,0,0)x_{1}^{2}+\sum\limits_{t=0}^{t=n+1}c_{\xi t}(x_{0},x_{1},0,0,0)x_{0}^{t}x_{1}^{n+2-t}+c_{l1}(x_{0},x_{1},0,0,0)x_{0}^{n+1}x_{1}.$$
Since each term of R.H.S of the equation contains $x_{1}$, but L.H.S does not contain $x_{1}$, 
we see an absurd situation. Therefore the above equation is not possible. \qed

\bibliographystyle{amsalpha}

\end{document}